\documentclass[12pt,english]{article}

\usepackage{amssymb}
\usepackage{polski}
\usepackage{ifthen}
\selecthyphenation{polish}
\mathchardef\ogon="012C%
\newcommand{\as}{a\kern-0.22em\lower.40ex\hbox{$_{\ogon}$}}
\newcommand{\As}{A\kern-0.22em\lower.40ex\hbox{$_{\ogon}$}}
\newcommand{\es}{e\kern-0.24em\lower.40ex\hbox{$_{\ogon}$}}
\newcommand{\Es}{E\kern-0.22em\lower.40ex\hbox{$_{\ogon}$}}
\makeatletter

\usepackage[a4paper, left=2.5cm, right=2.5cm, top=2.5cm, bottom=3.5cm, headsep=1.2cm]{geometry}
\newtheorem{theorem}{Theorem}[section]

\newtheorem{definition}[theorem]{Definition}
\newtheorem{example}[theorem]{Example}

\newtheorem{lemma}[theorem]{Lemma}

\newtheorem{proposition}[theorem]{Proposition}
\newtheorem{remark}[theorem]{Remark}

\newenvironment{Proof}[1][Proof]{\noindent\textbf{#1.} }{\ \rule{0.5em}{0.5em}}
\def\qed{\hbox to 0pt{}\hfill$\rlap{$\sqcap$}\sqcup$}
\makeatother
\usepackage{amssymb}
\usepackage{amsmath}
\usepackage{amsfonts}
\usepackage{color}
\usepackage{graphicx}
\usepackage{graphics}
\numberwithin{equation}{section}
\usepackage{babel}
\date{}
\title{On global properties of lower semicontinuous quadratically minorized functions }
\author{Monika Syga
	\thanks{Warsaw University of Technology, Faculty of Mathematics and Information Science, ul. Koszykowa 75,
		00--662 Warsaw, Poland,  M.Syga@mini.pw.edu.pl} 
}
\begin{document}
	\maketitle

\begin{abstract}%
We use the framework of a type of abstract convexity ($\Phi_{lsc}$-convexity) to investigate properties of   lower semicontinuous quadratically minorized functions in Hilbert spaces.  A new result, which states  that, for every local $\Phi_{lsc}$-subgradient there exists a global one is proved and plays a crucial role in our considerations. We deliver conditions for abstract subdifferentiability ($\Phi_{lsc}$-subdifferentiability) of  locally $C^{1,1}$ functions, twice continuously differentiable functions, prox-regular functions and paraconvex functions.  As an application we establish a new sufficient and necessary condition for minimax equality for $\Phi_{lsc}$-convex functions. This  new condition is expressed in therms of $\Phi_{lsc}$-subdifferential.


\medskip{}

\textbf{Mathematics Subject Classification (2000)32F17; 49J52; 49K27; 49K35; 52A01}

\end{abstract}

\section{Introduction}

In the present  paper we use the tools of so called $\Phi$-convexity to
study the properties of  lower semicontinuous quadratically minorized functions (l.s.q.m for short).
Lower semicontinuous quadratically minorized functions  appear  in many context, e.g. in approximation theory \cite{Attouch-Aze} or in image processing \cite{Angulo}.
An important examples of l.s.q.m functions are prox-bounded functions \cite{poli-rock96} and paraconvex functions \cite{Rolewiczpara}, known also under the name weakly convex functions  \cite{Vial} and semiconvex functions \cite{cannarsa}.
The $\Phi$-convexity theory (\cite{rolewicz},\cite{rubbook}), a type of abstract convexity, provides a suitable framework which allows to treat, in a unified way, different classes of functions, mainly  in context of global optimization problems (e.g. \cite{ioffe-rub}, \cite{rub-hua}).  $\Phi$-convex functions are defined as pointwise suprema of functions from a given class $\Phi$. Such an approach
to abstract convexity generalizes the classical fact that each 
proper  lower semicontinuous 
convex function is the upper envelope of a certain set of affine functions.  

The present paper is devoted to the study of $\Phi_{lsc}$-subdifferential  and the  $\Phi_{lsc}$-sub-\\differentiability of a given $\Phi_{lsc}$-convex function, where the class $\Phi_{lsc}$ is defined as
\begin{equation}
\tag{*}
\label{lscset}
\Phi_{lsc}:= \{\varphi : X \rightarrow \mathbb{R}, \ \varphi(x)=-a\|x\|^2+ \left\langle v,x\right\rangle+c, \ \ x\in X,\  v\in X^*, \ a\geq 0, \ c\in \mathbb{R} \},
\end{equation}
with $X$ a Hilbert space and $X^*$ its topological dual.

It is well known feature of convex functions that every local subgradient (in the sens of the convex analysis) is also a global one.  In \cite{rol-glob} such property for a general $\Phi$-convex functions is called the {\em globalization property}, and studied for a number of classes $\Phi$ therein. 
We show that the class of $\Phi_{lsc}$-convex functions possess the globalization property  i.e.  the existence   of a local $\Phi_{lsc}$-subgradient of a function $f$ at a given point  implies the existence of a global $\Phi_{lsc}$-subgradient (Proposition \ref{glob-prox}).   This result is crucial and is the main tool used in the paper.
With help of the globalization property we establish the $\Phi_{lsc}$-subdifferentiability of  many important classes of functions, such  as  locally $C^{1,1}$ functions, twice continuously differentiable functions and paraconvex functions. Moreover we show that in the class of  prox-regular functions (which appears frequently in variational analysis and covers functions such as strongly amenable, lower-$C^2$ and primal-lower-nice (pln) \cite{RockWets98}) the $\Phi_{lsc}$-subdifferentiability is assured. We also discus the connections between the $\Phi_{lsc}$-subdifferential and number of already known generalized subdifferentials e.g. proximal subdifferential, Clarke and Dini subdifferential.

The main contribution of the paper lies in the use of  $\Phi_{lsc}$-convexity, which enables us to treat many well known classes of functions in the unifying way and thus establish some known and new results.  

To show an application of the  $\Phi_{lsc}$-subdifferential we deliver conditions, involving this subdifferential, for the minimax equality
$$
\sup_{y\in Y} \inf_{x\in X} a(x,y)=\inf_{x\in X} \sup_{y\in Y} a(x,y),
$$
where $X,Y$ are nonempty sets and $a:X\times Y\rightarrow\hat{\mathbb{R}}:=\mathbb{R}\cup\{\pm\infty\}$ is lower semicontinuous and quadratically minorized  as a function of $x$ and concave (in the classical sense) as a function of $y$. An exhaustive survey of minimax theorems is given e.g. in \cite{simons}.
According to our knowledge,  in the literature, there is no minimax theorems addressing directly l.s.q.m functions.
This class of functions appear frequently in optimization problems and the respective minimax theorems can be used  to provide weak duality theorems for optimization problems involving l.s.q.m. functions. It is an important issue  to provide, as weak as possible, conditions under which the minimax equality for such functions holds. To study these conditions we use the general minimax theorem for $\Phi$-convex functions (\cite{Syga2018}).

If, for every $y\in Y$  the function $a(\cdot,y)$ is $\Phi_{lsc}$-convex then, a sufficient and  necessary
condition for $a(\cdot,\cdot)$ to satisfy the minimax equality is so called intersection property
introduced in \cite{bed-syg} and investigated in \cite{bed-syg-joca, syga, Syga2018}. 
\begin{definition}
	\label{ip}
	Let $\varphi_1,\varphi_2:X\rightarrow \mathbb{R}$ be any functions from the set $\Phi_{lsc}$ \eqref{lscset} and $\alpha\in \mathbb{R}$. We say that 
	{\em the intersection property} holds for $\varphi_1$ and $\varphi_2$ on $X$ at the level $\alpha$ if and only if
	\begin{equation}
	\label{empty}
	\tag{**}
	[\varphi_1<\alpha]\cap[\varphi_2<\alpha]=\emptyset,
	\end{equation}
	where $[\varphi<\alpha]:=\{x\in X\ :\ \varphi(x)<\alpha\}$ is the strict lower level set 
	of function $\varphi:X\rightarrow \mathbb{R}$.
\end{definition}
\vskip 0.15 true in
Let us note that the intersection property is   expressed via an algebraic condition \eqref{empty} which  is not easy to check.  For instance, sufficient conditions for the intersection property in the convex case  are given in Theorem 5.2 of \cite{bed-syg-joca} (see also Theorem 4.5 of \cite{bed-syg-joca} for some relationships between the intersection property and subdifferentials). At this point it is worth observing that conditions similar to \eqref{empty} appear in other context, e.g. in so called S-lemma \cite{s-lemma}.

In this paper  we investigate the relationships between the intersection property stated above and $\Phi_{lsc}$-subgradients of functions $a(\cdot,y)$ (Proposition  \ref{prop1n}). We  introduce a new condition, called {\em zero subgradient condition},  which is a sufficient condition for the intersection property (Definition \ref{ip}).

The organization of the paper is as follows. In the next section we present basic notions, definitions and  properties of the class of $\Phi_{lsc}$-convex functions and its important subclasses (such as paraconvex functions). Section 3 is the main section of the paper and is devoted to study  $\Phi_{lsc}$- subdifferential. In Proposition \ref{glob-prox} we show that if there exits a local $\Phi_{lsc}$-subgradient then there exists a global one. Then we make use of this result and establish the $\Phi_{lsc}$-subdifferentiability of paraconvex functions (Proposition \ref{par}),    locally $C^{1,1}$ functions (Proposition \ref{gat}),   twice continuously differentiable functions (Proposition \ref{fen}) and prox-regular functions (Proposition \ref{prox-diff}). Moreover, we establish the relationship between the $\Phi_{lsc}$-subdifferential and Dini and Clarke subdifferentials (Proposition \ref{dini-clark}).

In  section 4 we introduce the zero subgradient condition and we show that this condition is a sufficient condition for the intersection property (Proposition \ref{prop1n}), we also prove that this condition in slightly modified from is a necessary condition for the intersection property (Proposition  \ref{nes}). The last section is devoted to minimax theorems for $\Phi_{lsc}$-convex functions (Theorem \ref{new_min_max_1}, \ref{new_min_max_2}, \ref{ness_min_max} and  \ref{new_min_max_3}).
\section{$\Phi_{lsc}$-convexity - the unifying framework}
In the present section we use the tools of $\Phi$-convexity to deal with very broad class of lower semicontinuous quadratically minorized functions defined on Hilbert space. We start by recalling basic notions and definitions. Throughout the paper  $X$ is  a Hilbert space with inner product $\langle \cdot,  \cdot\rangle :X\times X\rightarrow \mathbb{R}$ and the respective norm $\|\cdot\|:X\rightarrow \mathbb{R}$. An open ball centred at $x\in X$ is denoted by $B(\delta,x):=\{y\in X \ :\ \|y-x\|<\delta \}$, $\delta>0$.

As stated in the Introduction, the class  $\Phi_{lsc}$ is defined as
$$
\Phi_{lsc}:= \{\varphi : X \rightarrow \mathbb{R}, \ \varphi(x)=-a\|x\|^2+ \left\langle v,x\right\rangle+c, \ \ x\in X,\  v\in X^*, \ a\geq 0, \ c\in \mathbb{R} \}.
$$
For any $f:X\rightarrow\hat{\mathbb{R}}:=\mathbb{R}\cup \{-\infty \}\cup \{ +\infty \}$  the set
$$
\text{supp}(f):=\{\varphi\in \Phi_{lsc}\ :\ \varphi\le f\}
$$
is called the {\em support} of $f$, where $\varphi\le f$ is defined as $\varphi(x)\le f(x)$ for all $x\in X$.  Whenever we say that $f$ is minorized by a quadratic function, it means that there exists $\bar{\varphi}\in \Phi_{lsc}$ such that $f\geq \bar{\varphi}$, i.e. the set $\text{supp}(f)$ is nonempty.
\begin{definition}(\cite{dolecki-k, rolewicz, rubbook})
	\label{convf}
	A function $f:X\rightarrow
	\hat{\mathbb{R}}$ is called {\em $\Phi_{lsc}$-convex} on X if
	$$
	f(x)=\sup\{\varphi(x)\ :\ \varphi\in\textnormal{supp}(f)\}\ \ \forall\ x\in X.
	$$
\end{definition}
By convention, if $f\equiv-\infty$ then $\textnormal{supp}(f)=\emptyset$.
In this paper we  limit our attention to functions $f:X\rightarrow\bar{\mathbb{R}}:=\mathbb{R}\cup\{+\infty\}$ such that $\textnormal{supp}(f)\neq \emptyset$.
We say that a function $f:X\rightarrow\bar{\mathbb{R}}$ is proper if  $\textnormal{supp}(f)\neq \emptyset$  and the effective domain of $f$ is nonempty, i.e.
$$
\text{dom}(f):=\{x\in X \ : \ f(x)<+\infty \}\neq \emptyset.
$$

It is a well known result that a function defined on Hilbert space is lower semicontinuous if and only if is supremum of continuous functions.  
In the following theorem we recall a characterization of $\Phi_{lsc}$-convex functions.

\begin{proposition}
	\label{phi-convex}
	Let $f:X\rightarrow\bar{\mathbb{R}}$ be a proper function. $f$ is  $\Phi_{lsc}$-convex on $X$ if and only if $f$ is lower semicontinuous on $X$ and minorized by a function from the class $\Phi_{lsc}$.
\end{proposition}

\begin{Proof}
	$\Rightarrow$ Let $f$ be  $\Phi_{lsc}$-convex function. Since the class $\Phi_{lsc}$ consists of continuous functions we get that $f$ is lower semicontinuous. By contradiction, assume that $f$ is not minorized by a function from the class $\Phi_{lsc}$, this means that $\text{supp}(f)=\emptyset$. It is a contradiction with the assumption that $f$ is proper.
	
	$\Leftarrow $ Let $f$ be lower semicontinuous and minorized by a function from the class $\Phi_{lsc}$ i.e.  $\textnormal{supp}(f) \neq \emptyset$, from \cite{rubbook}, Example 6.2 we get that $f$ is $\Phi_{lsc}$-convex. 
\end{Proof}

Now we discuss several important subclasses of $\Phi_{lsc}$-convex functions, which often appear in
applications.  We start with so called $\gamma$-paraconvex functions, 
which were first introduced in \cite{rolewicz1979}, considered e.g. in \cite{daniilidis} and appear in context of optimization in e.g. \cite{BBEM}, Let $\gamma $ be a positive number.
\begin{definition}
	A function $f:X\rightarrow\bar{\mathbb{R}}$ is called {\em $\gamma$-paraconvex} on $X$ if there exists $C>0$ such that for all $x,y\in X$  and $t\in [0,1]$ the following inequality holds
	\begin{equation}
	\label{para}
	f(tx+(1-t)y)\leq tf(x)+(1-t)f(y)+C \|x-y\|^{\gamma}.
	\end{equation}
\end{definition}
Usually $2$-paraconvex functions  are called paraconvex.
The strong $\gamma$-paraconvexity was defined in \cite{rolewicz2000}. 
\begin{definition}
	A function $f:X\rightarrow\bar{\mathbb{R}}$ is {\em strongly $\gamma$-paraconvex} on $X$ if there exists $C>0$ such that for all $x,y\in X$  and $t\in [0,1]$ the following inequality holds
	\begin{equation}
	\label{spara}
	f(tx+(1-t)y)\leq tf(x)+(1-t)f(y)+C\min\{t,1-t \} \|x-y\|^{\gamma}.
	\end{equation}
\end{definition}
It is obvious that if a function is strongly $\gamma$-paraconvex then it is $\gamma$-paraconvex. It was shown in \cite{rolewicz1979} that for $\gamma \in (1,2]$ the strong  $\gamma$-paraconvexity is equivalent to $\gamma$-paraconvexity. In context of our consideration we focus on $2$-paraconvexity, which throughout the paper will be called paraconvexity.

The class of paraconvex functions coincides with other classes investigated in the literature. Now we will discuss characterizations of paraconvex functions.
We start, with weakly convex functions, which were first introduced in \cite{Vial}, they appear in context of global optimization in \cite{Wu2007}, \cite{Wu2010} and in context of approximation theory in \cite{Attouch-Aze}.
\begin{definition}
	A function $f:X\rightarrow\bar{\mathbb{R}}$ is {\em weakly convex} on $X$ if there exists $c>0$ such that the function $f(x) +c\|x\|^2$ is convex. 
\end{definition}

In \cite{rolewicz} the following proposition was shown.

\begin{proposition}[\cite{rolewicz}, Proposition 5.2.11]
	\label{pr-wp}
	Let $f:X\rightarrow\bar{\mathbb{R}}$ be a proper function. $f$ is  weakly convex on $X$ if and only if  $f$ is paraconvex on $X$.
\end{proposition}
In some papers weakly convex functions are also called semiconvex (see i.e. \cite{cannarsa}, Proposition 1.1.3). 
Let us note that in  Definition 4 and Definition 5 the function  $f$ is an arbitrary  whereas  in \cite{rolewicz} paraconvexity is defined for continuous functions only. 
Hence, the classes of  weakly convex, semiconvex, paraconvex and strongly paraconvex functions coincide. In the sequel, functions from any of these classes will be refereed to as  paraconvex.

The following proposition shows that the class of paraconvex functions (strongly paraconvex, weakly convex, semiconvex) is a subclass of $\Phi_{lsc}$-convex functions. 
\begin{proposition}
	Let $f:X\rightarrow\bar{\mathbb{R}}$ be a proper lower semicontinuous function. If $f$ is paraconvex on $X$ then $f$ is $\Phi_{lsc}$-convex on $X$. 
\end{proposition}
\begin{Proof}
	By the paraconvexity of $f$, there exists $c>0$ such that $f+c\|\cdot\|^2$ is convex on $X$. Since  $f+\|\cdot\|^2$ is lower semicontinuous, by  Proposition 3.1, \cite{Ekeland},   $f+\|\cdot\|^2$ can be represented as pointwise suprema of affine functions  i.e
	$$
	f(x)+c\|x\|^2=\sup\{\langle v,x\rangle+b \ : \langle v,x\rangle+b \leq	f(x)+c\|x\|^2 \ v\in X^*, b\in \mathbb{R}\} \ \ \forall\ x\in X.
	$$
	Consequently,
	$$
	f(x)=\sup\{-c\|x\|^2+\langle v,x\rangle+b
	\ : -c\|x\|^2+\langle v,x\rangle+b \leq	f(x),  \ v\in X^*, b\in \mathbb{R}\}\ \ \forall\ x\in X,
	$$
	i.e. $f$ is $\Phi_{lsc}$-convex.
\end{Proof}

In the next section we will show examples of functions which are $\Phi_{lsc}$-convex but not paraconvex. 

Following \cite{RockWets98}, we say that a proper lower semi-continuous function $f:X\rightarrow\bar{\mathbb{R}}$ is {\em prox-bounded} if there exists a polynomial  $q$ od degree two or less such that $f\geq q$ (see \cite{RockWets98}, Exercise 1.24). By Proposition \ref{phi-convex}, in Hilbert spaces the class of prox-bounded functions coincides with the class of $\Phi_{lsc}$-convex functions. The detailed study of prox-bounded functions can be found in \cite{RockWets98} for finite dimensional spaces and in \cite{BernardThibault2005} for Hilbert spaces.

Let us note that the set of all $\Phi_{lsc}$-convex functions defined on normed space $X$ contains all proper lower semicontinuous and convex (in the classical sense) functions defined on $X$.

In further considerations we will focus on the class of $\Phi_{lsc}$-convex functions and its important subclass which, in order to avoid repetition, will be called the class of  paraconvex functions.  By above consideration, all the results in the next sections, which are proved for $\Phi_{lsc}$-convex functions,  are also true for prox-bounded functions and  paraconvex functions which coincides with the classes of strongly paraconvex, weakly convex and semiconvex functions. 

\section{$\Phi_{lsc}$-subdifferential}
We  cast the concept of $\varepsilon$-$\Phi_{lsc}$-subdifferential into our framework of  $\Phi_{lsc}$-convexity. The $\varepsilon$-$\Phi$-subdifferential, for an arbitrary class $\Phi$, was defined in \cite{ioffe-rub}, which is a direct adaptation of the classical  definition of $\varepsilon$-subdifferential for a convex function. Let $\mathbb{R}_{+}$ be a set of all nonnegative numbers.
\begin{definition}
	Let $f:X\rightarrow\bar{\mathbb{R}}$ be a proper function and let $\varepsilon\ge 0$.  An element  $(a,v)\in \mathbb{R}_+\times X^*$ is called a {\em $\varepsilon$-$\Phi_{lsc}$-subgradient} of $f$ at $\bar{x}\in \text{dom}(f)$, if the following inequality holds
	\begin{equation}
	\label{phinew}
	f(x)-f(\bar{x}) \geq 	\langle v, x-\bar{x}\rangle -a\| x\|^2+a\|\bar{x}\|^2-\varepsilon,  \  \ \ \ \ \ \forall \ x\in X .
	\end{equation}
	The set of all  $\varepsilon$-$\Phi_{lsc}$-subgradients of $f$ at  $\bar{x}$ is denoted as $\partial_{lsc}^{\varepsilon} f(\bar{x})$, if $\varepsilon=0$, then we write $\partial_{lsc}f(\bar{x})$. If $(0,0)\in \partial_{lsc}^{\varepsilon} f(\bar{x})$, then we simply write $0 \in  \partial_{lsc}^{\varepsilon} f(\bar{x})$. Function $f$ is $\Phi_{lsc}$-subdifferentiable at $\bar{x}$ if  $ \partial_{lsc} f(\bar{x})\neq \emptyset$.
\end{definition}

\begin{remark}
	\label{rem1}
	\begin{description}
		\item{a)} It is easy to show that for a proper function $f$ the set $\partial_{lsc}^{\varepsilon}  f(\bar{x})$ is convex for all $\varepsilon\geq 0$ and $\bar{x}\in \text{dom}(f)$.
		\item{b)} For a proper $\Phi_{lsc}$-convex function $f$ and for every $\varepsilon>0$ and $\bar{x}\in \text{dom}f$ the set $\partial_{lsc}^{\varepsilon} f(\bar{x})$ is nonempty. Indeed, let $\varepsilon>0$. By $\Phi_{lsc}$-convexity of $f$ for an arbitrary $\bar{x}\in \text{dom}f$ we have
		$$
		f(\bar{x})=\sup\{\varphi(\bar{x})\ :\ \varphi\in\textnormal{supp}(f)\},
		$$
		hence,  there exists $\bar{\varphi}\in \textnormal{supp}(f) $ such that 
		$$
		\bar{\varphi}(\bar{x})> f(\bar{x})-\varepsilon.
		$$
		Consequently,
		$$
		f(x)-f(\bar{x}) \geq \bar{\varphi}(x)-\bar{\varphi}(\bar{x}) -\varepsilon\ \ \ \ \ \forall x\in X.
		$$
		Since $\bar{\varphi} \in \Phi_{lsc}$, there exist $a\geq0$, $v\in X^*$ and $c\in \mathbb{R}$ such that $\bar{\varphi}(\cdot)=-a\|\cdot\|^2+ \left\langle v,\cdot\right\rangle+c$. Hence, the above inequality is equivalent to
		$$
		f(x)-f(\bar{x}) \geq 	\langle v, x-\bar{x}\rangle -a\| x\|^2+a\|\bar{x}\|^2-\varepsilon,  \  \ \ \ \ \ \forall \ x\in X,
		$$
		i.e. $(a,v)\in\partial_{lsc}^{\varepsilon} f(\bar{x})$.
		\item{c)} If for a proper  $\Phi_{lsc}$-convex function $f$ there exists $\bar{x}\in \text{dom}(f)$ such that the supremum is attained i.e. there exists $\bar{\varphi}\in \text{supp}(f)$ such that $f(\bar{x})=\bar{\varphi}(\bar{x})$, then $ \partial_{lsc}f(\bar{x}) \neq \emptyset$.
	\end{description} 
\end{remark}
\vskip 0.15 true in
Now we recall the definition of a local $\Phi$-subgradient. This concept for general class $\Phi$ was discussed in details e.g. in \cite{rol-glob}.
\begin{definition}
	Let $f:X\rightarrow\bar{\mathbb{R}}$ be a proper function. An element  $(a,v)\in \mathbb{R}_+\times X^*$ is called a {\em local $\Phi_{lsc}$-subgradient} of $f$ at $\bar{x}\in \text{dom}(f)$, if there exists $\delta >0$ such that, the following inequality holds
	\begin{equation}
	\label{phinew1}
	f(x)-f(\bar{x}) \geq 	\langle v, x-\bar{x}\rangle -a\| x\|^2+a\|\bar{x}\|^2,  \  \ \ \ \ \ \forall \ x\in B(\delta,\bar{x}) .
	\end{equation}
	The set of all  local $\Phi_{lsc}$-subgradients of $f$ at  $\bar{x}$ is denoted by $\partial_{lsc}^{\text{loc}} f(\bar{x})$.
\end{definition}
For a given class $\Phi$, the fact that a function $f$ has a local $\Phi$-subgradient at a point $\bar{x}$ dose not imply the existence of a global $\Phi$-subgradient. However, there are classes of $\Phi$-convex functions with this property (for a number of examples see  \cite{rol-glob}). We now show that for a class of $\Phi_{lsc}$-convex functions the existence of  a local $\Phi_{lsc}$-subgradient at a point $\bar{x}$ indicate the existence of a global $\Phi_{lsc}$-subgradient. To this aim we first prove the following lemma (the idea of the prove is based on that of  Lemma 3.7 in \cite{BBEM}).
\begin{lemma}(\cite{BBEM}, Lemma 3.7)
	\label{lembbem}
	Let $f:X\rightarrow \bar{\mathbb{R}}$ be a proper $\Phi_{lsc}$-convex function and $\bar{x}\in \text{dom}(f)$. If  there exists $\rho\geq$ and $\delta>0$ such that
	\begin{equation}
	\label{lemm1}
	f(x)\geq f(\bar{x})-\rho\|x-\bar{x}\|^2 \ \ \ \ \forall \ \ x\in B(\delta, \bar{x})
	\end{equation}
	then there exists $\bar{\rho}\geq 0$ such that
	\begin{equation}
	\label{p-zero}
	f(x)\geq f(\bar{x})-\bar{\rho}\|x-\bar{x}\|^2 , \ \ \ \ \forall \ \ x\in X.
	\end{equation}
\end{lemma}
\begin{Proof}
	Let $\rho\geq$ and $\delta>0$ be such that the inequality \eqref{lemm1} holds.
	It is easy to see that,	by the $\Phi_{lsc}$-convexity of the function $f$, there exists $a\geq 0$ and $c\in \mathbb{R}$ such that
	\begin{equation}
	\label{eq1a}
	f(x)\geq -a\|x\|^2 +c \ \ \ \ \ \ \ \ \ \forall \ \ x\in X.
	\end{equation}
	For all $x\in X$ we have 
	$$
	f(x)	\begin{array}[t]{l}\geq -a\|x\|^2 +c\\
	=-a\|x-\bar{x}\|^2-a\|\bar{x}\|^2+2a\langle \bar{x}-x,\bar{x}\rangle+c\\
	\geq -a\|x-\bar{x}\|^2-a\|\bar{x}\|^2-2a\| x-\bar{x}\| \|\bar{x}\|+c.
	\end{array}
	$$
	Hence, for all $x\notin B(\delta, \bar{x})$, i.e. $\|x-\bar{x}\|\geq \delta$ we have the following inequality
	\begin{equation}
	\label{eq1}
	f(x)	\begin{array}[t]{l} \geq -a\|\bar{x}\|^2 +c -a\|x-\bar{x}\|^2\left( 1+\frac{2\|\bar{x}\|}{\| x-\bar{x}\| }\right)\\
	\geq 	-a\|\bar{x}\|^2 +c -a\|x-\bar{x}\|^2\left( 1+\frac{2\|\bar{x}\|}{\delta}\right).
	\end{array}
	\end{equation}
	Let 
	$$
	\bar{a}:= \frac{f(\bar{x}) +a\|\bar{x}\|^2 -c}{\delta^2}+a\left( 1+\frac{2\|\bar{x}\|}{\delta}\right)
	$$
	By \eqref{eq1a}, we have $	f(\bar{x})\geq -a\|\bar{x}\|^2 +c$, hence $\bar{a}\geq 0$. By the definition of $\bar{a}$, we have
	$$
	-a\|\bar{x}\|^2 +c -a\|x-\bar{x}\|^2\left( 1+\frac{2\|\bar{x}\|}{\delta}\right)\geq f(\bar{x})-\bar{a}\|x-\bar{x}\|^2,
	$$
	for all $x\notin B(\delta, \bar{x})$.
	From \eqref{eq1}, the inequality
	$$
	f(x)\geq f(\bar{x})-\bar{a}\|x-\bar{x}\|^2,
	$$
	holds for all $x\notin B(\delta, \bar{x})$. Hence, inequality \eqref{p-zero} holds with $\bar{\rho}=\max\{ \bar{a},\rho  \}$.
\end{Proof}

Now we use Lemma \ref {lembbem} to prove that the existence of  a local $\Phi_{lsc}$-subgradient at a point $\bar{x}$ implies the existence of a global $\Phi_{lsc}$-subgradient.
\begin{proposition}
	\label{glob-prox}
	Let $f:X\rightarrow \bar{\mathbb{R}}$ be a proper $\Phi_{lsc}$-convex function and $\bar{x}\in \text{dom}(f)$. If $(a,v)\in\partial_{lsc}^{loc}f(\bar{x})$ then there exists $\bar{a}\geq 0$  such that $(\bar{a},v-2a\bar{x}+2\bar{a}\bar{x} )\in\partial_{lsc}f(\bar{x})$.
\end{proposition}
\begin{Proof}
	Let $(a,v)\in\partial_{lsc}^{loc}f(\bar{x})$, hence there exists $\delta >0$  such that
	$$
	f(x)-f(\bar{x}) \geq 	\langle v, x-\bar{x}\rangle -a\| x\|^2+a\|\bar{x}\|^2,  \  \ \ \ \ \ \forall \ x\in B(\delta,\bar{x}) .
	$$
	We have
	\begin{equation}
	\label{new1}
	f(x)-f(\bar{x}) \geq 	\langle v, x-\bar{x}\rangle -a\| x\|^2+a\|\bar{x}\|^2+2a\langle x-\bar{x}, \bar{x}\rangle -2a\langle x-\bar{x}, \bar{x}\rangle 
	\end{equation}
	$$
	= \langle v-2a\bar{x}, x-\bar{x}\rangle -a\|x-\bar{x}\|^2,
	$$
	for all $x\in B(\delta,\bar{x})$.  
	Let $h(\cdot):=f(\cdot)-\langle v-2a\bar{x}, \cdot-\bar{x}\rangle$. The function $h$ is $\Phi_{lsc}$-convex.  Indeed, by the $\Phi_{lsc}$-convexity of $f$,  there exists $\varphi\in \Phi_{lsc}$ such that
	$$
	f(x)\geq \varphi(x) \ \ \ \ \ \ \forall \ \  x\in X.
	$$
	Consequently,
	$$
	h(x)\geq \varphi(x)-\langle v-2a\bar{x}, x- \bar{x}\rangle \ \ \ \ \ \ \forall \ \  x\in X,
	$$
	and, the function $\bar{\varphi}(\cdot):=  \varphi(\cdot)-\langle v-2a\bar{x}, \cdot- \bar{x}\rangle $, belongs to  $\Phi_{lsc}$,  $\bar{\varphi}\in\Phi_{lsc}$. Furthermore,  lower semicontinuity of $f$ implies  lower semicontinuity of $h$, since $\langle v-2a\bar{x}, \cdot- \bar{x}\rangle$ is a continuous function.
	Moreover, for the function $h$, the inequality \eqref{new1} takes the form
	$$
	h(x)\geq h(\bar{x}) -a\|x-\bar{x}\|^2 , \ \ \ \ \forall \ \ x\in B(\delta,\bar{x}).
	$$
	Applying Lemma \ref{lembbem} to the function $h$, we get that there exists $\bar{a}\geq0$ such that
	$$
	h(x)\geq h(\bar{x}) -\bar{a}\|x-\bar{x}\|^2 , \ \ \ \ \forall \ \ x\in X.
	$$
	By the definition of $h$,  we get 
	$$
	f(x)\geq f(\bar{x})+\langle v-2a\bar{x}, x- \bar{x}\rangle  -\bar{a}\|x-\bar{x}\|^2 , \ \ \ \ \forall \ \ x\in X,
	$$
	which is equivalent to 
	$$
	f(x)\geq f(\bar{x}) + \langle v-2a\bar{x}+2\bar{a}\bar{x}, x- \bar{x}\rangle  -\bar{a}\| x\|^2+\bar{a}\|\bar{x}\|^2, \ \ \ \ \forall \ \ x\in X.
	$$
	i.e. $(\bar{a},v-2a\bar{x}+2\bar{a}\bar{x} )\in\partial_{lsc}f(\bar{x})$.	\end{Proof}

\subsection{$\Phi_{lsc}$-subdifferentiability }
Let us note, that the nonemptiness of  set $\partial_{lsc}  f(x)$ on the domain of $f$ is not ensured even for differentiable functions. For example, the function $f:\mathbb{R}\rightarrow\mathbb{R}$, $f(x)=-|x|^\frac{3}{2}$ is $\Phi_{lsc}$-convex on $X$ and differentiable, but at the point $\bar{x}=0$ has no $\Phi_{lsc}$-subgradient. 
It turns out that  conditions ensuring the nonemptiness of $\Phi_{lsc}$-subdifferential can easily be formulated for paraconvex functions. The following proposition gives a simple criterion to distinguish  $\Phi_{lsc}$-convex functions which are not paraconvex.  

\begin{proposition}
	\label{par}
	Let $f:X\rightarrow\bar{\mathbb{R}}$ be a proper l.s.c. paraconvex function, then for every $x\in\textnormal{int}\,\textnormal{dom} (f)$ the set $\partial_{lsc}  f(x)$ is nonempty.
\end{proposition}
\begin{Proof}
	By the paraconvexity of $f$, there exists $c>0$ such that the function
	$$
	g(\cdot):=f(\cdot)+c\|\cdot\|^2
	$$
	is convex. By the lower semicontinuity of $f$ we get that $g$ is lower semicontinuous. It is well known result from convex analysis that $g$ have the classical subgradient at every point of $\textnormal{int}\,\text{dom} (g)$ (see e.g. Theorem 2.4.9 of \cite{zalinescu2002}), i.e for every $\bar{x}\in\textnormal{int}\,\text{dom} (g)$ there exists $v\in X^*$ such that
	$$
	g(x)-g(\bar{x})\geq \langle v, x-\bar{x} \rangle \ \ \ \ \ \ \ \ \forall x\in X.
	$$
	By the definition of $g$, the above inequality is equivalent to
	$$
	f(x)-f(\bar{x})\geq \langle v, x-\bar{x} \rangle - c\|x\|^2 +c\|\bar{x}\|^2\ \ \ \ \ \ \ \ \forall x\in X,
	$$
	this means that $(c,v)\in \partial_{lsc}  f(\bar{x})$ and since $\text{dom}(f)=\text{dom}(g)$, the proof is completed.
\end{Proof}

Let us note that in the above Proposition it is possible to replace the $\textnormal{int}\,\textnormal{dom} (f)$ by  the so called quasi-relative interior  $\textnormal{qri}\,\textnormal{dom} (f)$ see i.e Corollary 9 in \cite{zalinescu2015}.
\vskip 0.15 true in
Taking into account Proposition \ref{par}, it is easy to see that, e.g. functions $f_1(x)=-|x|^\frac{3}{2}$, $f_2(x)=-|x|$ and $f_3(x)=||x|-1|$ are $\Phi_{lsc}$-convex, but  are not $\Phi_{lsc}$-subdifferentiable at the point $0$, which means that they are not paraconvex. On the other hand, if a function has $\Phi_{lsc}$-subgradient at every point of a open convex set  $U\subset X$, then $f$ is paraconvex on $U$.
This is the content of the following proposition.
\begin{proposition}
	\label{sub-para}
	Let $f:X\rightarrow\bar{\mathbb{R}}$ be a proper  function and $U\subset X$ be an open convex set such that $U\subset \textnormal{dom}(f)$. If there exists $a\geq 0$, such that  $(a,v_{\bar{x}})\in\partial_{lsc}f(\bar{x})$ for every $\bar{x}\in U$, then $f$ is paraconvex on $U$.
\end{proposition}
\begin{Proof}
	Let $\bar{x} \in U$  $(a,v_{\bar{x}})\in \partial_{lsc}f(\bar{x})$, then
	$$
	f(x)-f(\bar{x}) \geq 	\langle v_{\bar{x}}, x-\bar{x}\rangle -a\| x\|^2+a\|\bar{x}\|^2,  \  \ \ \ \ \ \forall \ x\in X.
	$$
	That inequality is equivalent to
	$$
	f(x)+a\| x\|^2-(f(\bar{x}) +a\|\bar{x}\|^2)\geq 	\langle v_{\bar{x}}, x-\bar{x}\rangle ,  \  \ \ \ \ \ \forall \ x\in X,
	$$
	which means that $v_{\bar{x}}$ is the classical (in the sense of the convex analysis) subgradient of the function $f+a\| \cdot\|^2$ at the point $\bar{x}$. Since  $\bar{x}\in U$ was chosen arbitrary, we have the conclusions that the function $f+a\| \cdot\|^2$ have a classical subgradient at every point of an open convex set $U$. Hence, $f+a\| \cdot\|^2$ is convex on $U$. Which means that $f$ is paraconvex on  $U$.
\end{Proof}
\vskip 0.15 true in
As we noted before, even differentiable $\Phi_{lsc}$-convex functions can have empty $\Phi_{lsc}$-subdifferential. Now we discuss classes of  differentiable (in the sense of  G\^{a}teaux or Fr\`{e}chet) functions which are also $\Phi_{lsc}$-subdifferentiable.  

By the $f'_{G}(x)$ we denote the G\^{a}teaux derivative of $f:X\rightarrow \mathbb{R}$ at a point $x\in  \text{dom}(f)$. We say that $f:X\rightarrow \mathbb{R}$  is locally $C^{1,1}$ around  $x\in \text{dom}(f)$ if  there exists $B(\delta,x)$ such that $f$ is G\^{a}teaux  differentiable for every $y\in  B(\delta, x)$ i.e  its G\^ateaux derivative  $f'_{G}(y)$ exists at every point $y\in  B(\delta,x)$,  and the G\^{a}teaux derivative $f'_{G}$ is Lipschitz continuous on $B(\delta,x)$. 
\begin{proposition}
	\label{gat}
	Let $f:X\rightarrow R$ be a proper $\Phi_{lsc}$-convex function. If $f$ is  $C^{1,1}$ around $\bar{x}\in \text{dom}(f)$, then there exists $\delta>0$ such that, for every $y\in B(\delta,\bar{x})$, the set $\partial_{lsc}f(y)$ is nonempty.
\end{proposition}
\begin{Proof}
	By assumption there exists $\delta>0$ such that, $f'_{G}$ is Lipschitz in $B(\delta,\bar{x})$ with the Lipschitz constant $\lambda >0$. Let $y\in B(\delta,\bar{x})$, for all $x\in B(\delta,\bar{x})$ and $t\in [0,1]$ we have
	
	$$
	| \langle   f'_{G}(tx +(1-t)y)- f'_{G}(y), x-y \rangle  |=\frac{1}{t}	| \langle  f'_{G}(t(x-y) +y)- f'_{G}(y), t(x-y) \rangle  |
	$$
	$$
	\leq t\lambda \|x-y\|^2.
	$$	
	Hence,  the following inequality holds
	$$
	f(x)-f(y)-\langle  f'_{G}(y), x-y \rangle = \int_{0}^{1} 	\langle  f'_{G}(tx +(1-t)y)- f'_{G}(y), x-y \rangle  \text{d}t
	\geq -\lambda \|x-y\|^2 \int_{0}^{1} t\text{d}t= -\frac{\lambda}{2} \|x-y\|^2.
	$$
	Which means
	$$
	f(x)-f(y)\geq \langle  f'_{G}(y), x-y \rangle -\frac{\lambda}{2} \|x-y\|^2 \ \ \ \forall \ x\in B(\delta,\bar{x}),
	$$
	equivalently 
	\begin{equation}
	\label{gat-1}
	f(x)-f(y)\geq \langle  f'_{G}(y)+\lambda y, x-y \rangle -\frac{\lambda}{2} \|x\|^2+\frac{\lambda}{2}\|y\|^2 \ \ \ \forall \ x\in B(\delta,\bar{x}),
	\end{equation}
	i.e. $(\frac{\lambda}{2},  f'_{G}(y)+\lambda y)\in \partial^{loc}_{lsc}f(y)$. By  Proposition \ref{glob-prox} there exists $a\geq 0$ such that $(a,f'_{G}(y)+2ay )\in  \partial_{lsc}f(y)$.
\end{Proof}
\begin{remark}
	The  inequality \eqref{gat-1}, means that the function $f(\cdot )+\frac{\lambda}{2}\|\cdot\|^2 $ is convex on $B({\delta},\bar{x})$. Hence, if $f$ is  $C^{1,1}$ around $\bar{x}\in \text{dom}(f)$, then there exists $\delta>0$ such that  $f$ is paraconvex on $B({\delta},\bar{x})$ with the constant $\frac{\lambda}{2}$. This results is well known in the literature, see i.e.  \cite{Drus}.
\end{remark}
\vskip 0.2 true in
By $f'_{F}(\bar{x})$ we denote the Fr\`{e}chet derivative of  $f:X\rightarrow \mathbb{R}$ at a point $x\in  \text{dom}(f)$.  We say that function $f$ is twice continuously differentiable on an open set $U\subset X$  if  $f$ is Fr\`{e}chet  differentiable at every point $x\in U$,  the  $f'_{F}$ is continuous function on $U$ and is also differentiable on $U$ with  the second derivative $f''_{F}$ being continuous function on $U$. The set of all twice continuously differentiable on $U$ functions is denoted as $C^2(U)$. For such functions we have the following Proposition.

\begin{proposition}
	\label{fen}
	Let $f:X\rightarrow R$ be a proper $\Phi_{lsc}$-convex function and $U$ be an open subset of $X$. If $f\in C^2(U)$, then for every $x\in U$, the set $\partial_{lsc}f(x)$ is nonempty.
\end{proposition}

\begin{Proof}
	Let $x\in U$.	By the fact  $f\in C^2(U)$,   there exists $B(\delta,x) $ such that $f$ admits a second-order Taylor expansion with a remainder i.e.
	\begin{equation}
	\label{weq1}
	f(y)=f(x)+\langle f'_{F}(x),y-x\rangle +\frac{1}{2}\langle f''_{F}(z)(y-x),y-x\rangle \ \ \ \forall \ \ y\in B(\delta,x),
	\end{equation}
	where $z$ is an element on the line segment between $x$ and $y$.
	
	Moreover, we have that there exists $\gamma >0$ such that 
	$$
	\langle f''_{F}(z)(y-x),y-x\rangle \geq -\| f''_{F}(z)\|\|y-x\|^2\geq -\gamma\|y-x\|^2 \ \ \ \forall \ \ \ z\in B(\delta,x).
	$$
	By the equality \eqref{weq1}, we have
	$$
	f(y)\geq f(x)+\langle f'_{F}(x),y-x\rangle -\frac{1}{2}\gamma\|y-x\|^2 \ \ \ \ \forall \ \ y\in B(\delta,x),
	$$
	equivalently
	\begin{equation}
	\label{f-1}
	f(y)\geq f(x)+\langle f'_{F}(x)+\gamma x,y-x\rangle -\frac{1}{2}\gamma\|y\|^2+\frac{1}{2}\gamma\|x\|^2 \ \ \ \ \forall \ \ y\in B(\delta,x),
	\end{equation}
	This means that $(\frac{1}{2}\gamma, f'_{F}(x)+\gamma x)\in \partial_{lsc}^{loc}f(x)$. By  Proposition \ref{glob-prox} there exists $a\geq 0$ such that $(a,f'_{F}(x)+2ax )\in  \partial_{lsc}f(x)$.
\end{Proof}
\begin{remark}
	The  inequality \eqref{f-1}, means that the function $f(\cdot )+\frac{\gamma}{2}\|\cdot\|^2 $ is convex on $B(\delta,x)$. Hence, if $f$ is  $C^{2}(U)$, then for every $x\in U$ there exists $\delta>0$ such that  $f$ is paraconvex on $B(\delta,x)$.
\end{remark}
\vskip 0.2 true in
By $\partial f(\bar{x})$ we denote {\em the limiting subdifferential} of  $f:X\rightarrow \mathbb{R}$ at a point $\bar{x}\in  \text{dom}(f)$ (see e.g. \cite{bor}).
We say that a proper function $f:X\rightarrow \mathbb{R}$ is {\em prox-regular} at a point $\bar{x}\in  \text{dom}(f)$ for $\bar{v}\in \partial f(\bar{x})$ if there exist $\rho>0$ and $\varepsilon >0$ such that for all $x,x'\in B(\bar{x},\varepsilon)$ with $|f(x)-f(\bar{x})|<\varepsilon$ and all $v\in B(\bar{v},\varepsilon)$ with $v\in\partial f (x)$ the following inequality holds
\begin{equation}
\label{prox-reg1}
f(x')\geq f(x) +	\langle v, x'-x\rangle -\frac{1}{2}\rho\|x'-x\|^2.
\end{equation}
If  $f$ is prox-regular for every $v\in \partial f(\bar{x})$, then $f$ is prox-regular at $\bar{x}$.  
For more details see \cite{poli-rock96} and \cite{BernardThibault2005} for the Hilbert space case. The class of prox-regular functions is an important class in variational analysis. This class contains strongly amenable, lower-$C^2$ and primal-lower-nice (pln) functions \cite{RockWets98}.  We prove that if function is prox-regular at a point $\bar{x}$ then is $\Phi_{lsc}$-subdifferentiable at $\bar{x}$. To this aim we use the concept of proximal subdifferential. A vector $v\in X^*$ is called a {\em proximal subgradient} of a function $f:X\rightarrow \bar{\mathbb{R}}$ at $\bar{x}\in \text{dom}(f)$ if there exists $\delta>0$ and $\rho\geq 0$ such that
\begin{equation}
\label{prox-in}
f(x)\geq f(\bar{x})+	\langle v, x-\bar{x}\rangle -\frac{1}{2}\rho\|x-\bar{x}\|^2 , \ \ \ \ \forall \ x\in B(\delta, \bar{x}).
\end{equation}

The set of all proximal subgradients of function $f$ at $\bar{x}$ is denoted as $\partial_{P}f(\bar{x})$ and called a proximal subdifferential. Proximal subgradients were studied e.g. in \cite{RockWets98}.

Now, we prove the equivalence between the $\Phi_{lsc}$-subdifferentiability and proximal subdifferentiability.
\begin{proposition}
	\label{prox-phi}
	Let $f:X\rightarrow \bar{\mathbb{R}}$ be a proper $\Phi_{lsc}$-convex function, $\bar{x}\in \text{dom}(f)$. Then
	$$
	\partial_{lsc}f(\bar{x})\neq \emptyset \ \Leftrightarrow \  	\partial_{P}f(\bar{x})\neq \emptyset.
	$$
\end{proposition}
\begin{Proof}
	Let $(a,w)\in \partial_{lsc}f(\bar{x})$, by the definition, for all $x\in X$ we have
	$$
	f(x)-f(\bar{x})\geq-a\|x\|^2+\langle w,x-\bar{x}\rangle +a\|\bar{x}\|^2.
	$$
	Let $\rho=2a$ and $v=w-\rho\bar{x}$, then $w=v+\rho\bar{x}$. For all $x\in X$ we have
	$$
	\begin{array}{l}
	f(x)-f(\bar{x})\geq \langle v+\rho\bar{x},x-\bar{x}\rangle -\frac{\rho}{2} \|x\|^2+\frac{\rho}{2} \|\bar{x}\|^2=\\
	\langle v,x-\bar{x}\rangle +\rho\langle\bar{x},x\rangle -\rho\langle\bar{x},\bar{x}\rangle-\frac{\rho}{2}\|x\|^2+\frac{\rho}{2} \|\bar{x}\|^2=\\
	\langle v,x-\bar{x}\rangle +\rho\langle\bar{x},x\rangle -\frac{\rho}{2}\|x\|^2-\frac{\rho}{2} \|\bar{x}\|^2=\\
	\langle v,x-\bar{x}\rangle-\frac{\rho}{2}\|x-\bar{x}\|^2,
	\end{array}
	$$
	i.e $v\in  \partial_{P}f(\bar{x})$ for every $\delta>0$.
	
	Let $v\in \partial_{P}f(\bar{x})$, by the definition there exists $\delta>0$ and $\rho\geq 0$ such that
	\begin{equation}
	\label{prox1}
	f(x)-f(\bar{x})\geq \langle v,x-\bar{x}\rangle -\frac{1}{2}\rho \|x-\bar{x}\|^2, \ \ \ \forall \ x\in B(\delta,\bar{x}).
	\end{equation}
	Equivalently,
	$$
	f(x)-f(\bar{x})\geq 	\langle v+\rho\bar{x},x-\bar{x}\rangle -\frac{1}{2}\rho \|x\|^2+\frac{1}{2}\rho \|\bar{x}\|^2\ \ \forall \ x\in B(\delta,\bar{x}),
	$$
	i.e. $(\frac{1}{2}\rho, v+\rho\bar{x})\in\partial_{lsc}^{loc}f(\bar{x})$. By  Proposition \ref{glob-prox}, there exists $\bar{a}\geq 0$  such that $(\bar{a},v+2 \bar{a}\bar{x})\in \partial_{lsc}f(\bar{x})$.
\end{Proof}

To  prove that prox-regular functions are $\Phi_{lsc}$-subdifferenttiable we use Proposition \ref{prox-phi} and the following theorem from \cite{BernardThibault2005}.
\begin{theorem}[\cite{BernardThibault2005}, Theorem 3.4]
	\label{thib}
	Let $f:X\rightarrow \bar{\mathbb{R}}$ be a proper l.s.c. function and $\bar{x}\in \text{dom}(f)$.  If  $f$ is prox-regular at  $\bar{x}$ for $\bar{v}\in \partial f(\bar{x})$ then $\bar{v}\in \partial_{P}f(\bar{x})$.
\end{theorem}
We have the following proposition
\begin{proposition}
	\label{prox-diff}
	Let $f:X\rightarrow \bar{\mathbb{R}}$ be a proper $\Phi_{lsc}$-convex function and $\bar{x}\in \text{dom}(f)$.  If  $f$ is prox-regular at  $\bar{x}$ for $\bar{v}\in \partial f(\bar{x})$  then there exists $a\geq 0$ such that $(\bar{a},\bar{v}+2 \bar{a}\bar{x})\in \partial_{lsc}f(\bar{x})$.
\end{proposition}
\begin{Proof}
	Follows immediately from Proposition \ref{prox-phi} and Theorem  \ref{thib}.
\end{Proof}

\subsection{Dini and Clarke subdifferential}
A natural question arises about the connection of  $\Phi_{lsc}$-subdifferential  with other known subdifferentials. This section is devoted to  Dini and Clarke  subdifferentials.
First, we recall some definitions. Let $f:X\rightarrow \bar{\mathbb{R}}$ be a proper function and $\bar{x}\in \text{dom}(f)$,. The {\em Clarke derivative} of  $f$ at a point $\bar{x}$ in a direction $h$ is 
$$
d^{C}f(\bar{x},h)=\limsup_{\substack{x\to \bar{x} \\t \downarrow 0}}{ \frac{f(x+th)-f(x)}{t}},
$$
where the upper limit is taken with respect to any sequence $\{t_n\}$ of positive number  tending to 0 and any sequence $\{x_n\}$ of elements belonging to $\text{dom}(f)$ with the limit $\bar{x}$. The following set is called the {\em Clarke subdifferential} of $f$ a the point $\bar{x}$
$$
\partial^{C}f(\bar{x})=\{v\in X^*\ :\ \langle v, h\rangle \leq d^{C}f(\bar{x},h) \ \forall h\in X  \}.
$$
The {\em Dini derivative} of  $f$ at a point $\bar{x}$ in a direction $h$ is 
$$
d^{D}f(\bar{x},h)=\liminf_{\substack{u\to h \\t \downarrow 0}}{ \frac{f(x+tu)-f(x)}{t}},
$$
where the lower limit is taken with respect to any sequence $\{t_n\}$ of positive number  tending to 0 and any sequence $\{u_n\}$ with the limit $h$. 
The following set is called the {\em Dini subdifferential} of $f$ a the point $\bar{x}$
$$
\partial^{D}f(\bar{x})=\{v\in X^*\ :\ \langle v, h\rangle \leq d^{D}f(\bar{x},h). \ \forall h\in X  \}
$$
Let us observe that
$$
d^{D}f(\bar{x},h) \leq d^{C}f(\bar{x},h),
$$
which implies that
$$
\partial^{D}f(\bar{x}) \subset \partial^{C}f(\bar{x}).
$$
The following Proposition holds
\begin{proposition}
	\label{dini-clark}
	Let $f:X\rightarrow \bar{\mathbb{R}}$ be a proper  $\Phi_{lsc}$-convex  function and $\bar{x}\in \text{dom}f$. Then
	$$
	\partial_{lsc}f(\bar{x})\neq \emptyset \ \Rightarrow  \ \partial^{D}f(\bar{x})\neq \emptyset \ \Rightarrow \ \partial^{C}f(\bar{x})\neq\emptyset.
	$$
\end{proposition}
\begin{Proof}
	We only need to show that $	\partial_{lsc}f(\bar{x})\neq \emptyset \ \Rightarrow  \ \partial^{D}f(\bar{x})\neq \emptyset$.	Let $(a,v)\in 	\partial_{lsc}f(\bar{x})$, then
	\begin{equation}
	\label{d1}
	f(x)-f(\bar{x}) \geq 	\langle v, x-\bar{x}\rangle -a\| x\|^2+a\|\bar{x}\|^2,  \  \ \ \ \ \ \forall \ x\in X .
	\end{equation}
	Let $h\in X$, $h_n\to h$ and $t_n\to 0$ be a sequence of positive numbers. Let $x=\bar{x}+t_nh_n$, the inequality \eqref{d1} takes the form
	$$
	f(\bar{x}+t_nh_n)-f(\bar{x}) \geq 	\langle v, t_nh_n\rangle +2a\langle \bar{x}, t_nh_n\rangle -a\|t_nh_n\|^2,
	$$
	hence 
	$$
	\frac{	f(\bar{x}+t_nh_n)-f(\bar{x})}{t_n} \geq 	\langle v+2a\bar{x}, h_n \rangle -at_n\|h_n\|^2.
	$$
	The last inequality means that the vector $v+2a\in \partial^{D}f(\bar{x}) $.
\end{Proof}

In general, for $\Phi_{lsc}$-convex function we cannot expect equivalences in Proposition \ref{dini-clark}. Consider the function $f:\mathbb{R}\rightarrow\mathbb{R}$, $f(x)=||x|-1|$, it is easy to see that $\partial_{lsc}f(0)=\emptyset$, on the other hand $\partial^{C}f(0)=[-1,1]$. 

It was shown by Jurani in \cite{jurani96} that the equivalences in Proposition \ref{dini-clark} hold if we limit ourselves to the class of paraconvex functions. 
Below we cite the result of Jurani, adopting the notation to our framework.
\begin{theorem}(\cite{jurani96}, Theorem 3.1)
	Let $f:X\rightarrow \bar{\mathbb{R}}$ be a  paraconvex function and $\bar{x}\in \text{dom}f$. 
	$$
	\partial_{lsc}f(\bar{x})\neq \emptyset \ \Leftrightarrow  \ \partial^{D}f(\bar{x})\neq \emptyset \ \Leftrightarrow \ \partial^{C}f(\bar{x})\neq\emptyset.
	$$
\end{theorem}

\section{Zero subgradient condition}
\label{grad} In this section we present  one of  possible applications of  the $\Phi_{lsc}$-subdifferential, i.e. conditions for minimax theorems for $ \Phi_{lsc}$-convex functions. To this aim for any two $\Phi_{lsc}$-convex functions, we introduce a condition, called {\em zero subgradient condition}.
This condition, expressed in terms of $\Phi_{lsc}$-subdifferentials, is  sufficient for the minimax equality to hold for $\Phi_{lsc}$-convex  functions (Theorem \ref{new_min_max_1}), and the modified version of this condition is also a necessary condition for the minimax equality (Theorem \ref{ness_min_max}).
Here we investigate properties of the zero subgradient condition  and  its relation  with the  intersection property defined in Definition \ref{ip}, which is a necessary and sufficient condition for the minimax equality for general $\Phi$-convex functions.
\begin{definition}
	\label{zsc}
	Let $f,g:X\rightarrow\bar{\mathbb{R}}$ be $\Phi_{lsc}$-convex functions, $x_{1}\in\text{dom}(f)$, $x_{2}\in\text{dom}(g)$  and $\varepsilon\geq 0$. We say that $f$ and $g$
	satisfy the {\em zero subgradient condition} at $(x_1,x_2)$, with $\varepsilon$ if 
	\begin{equation*}
	\label{main-con}
	0\in\mbox{co}(\partial^\varepsilon_{lsc}f(x_1)\cup\partial^\varepsilon_{lsc} g(x_2)),
	\end{equation*}
	where $\text{co}(\cdot)$ is a standard convex hull of a set. 
\end{definition}
For simplicity, if $f$ and $g$ satisfy the zero subgradient condition at ($x_1$, $x_2$) with $\varepsilon$ we will write that $f$ and $g$ satisfy the $ZS(\varepsilon,x_1,x_2)$ condition. If $x_1=x_2=\bar{x}$ we will write that $f$ and $g$ satisfy the $ZS(\varepsilon,\bar{x})$ condition.
\begin{remark}
	\begin{description}
		\item{a)}	Let us note that if, for a given function $f$, there exist $\bar{x}\in \text{dom}(f)$ and $\varepsilon\geq 0$, such that $0\in \partial^\varepsilon_{lsc}f(\bar{x})$, then $f$ and  every function $g$ defined on $X$ such that $\hat{x}\in \text{dom}(g)$, satisfy the $ZS(\varepsilon, \bar{x},\hat{x})$ condition, even if the set $\partial^\varepsilon_{lsc}g(\hat{x})$ is empty.
		
		\item{b)}	One can notice that the  $ZS(0,\bar{x})$ condition for $f$ and $g$ is similar to the formula for the  subdifferential of the function $\max\{f,g\}$ if $f(\bar{x})=g(\bar{x})$.
	\end{description}
\end{remark}
\vskip 0.1 true in
The following proposition shows that if two $\Phi_{lsc}$-convex functions satisfy the zero subgradient condition then,  one can find in their support sets two functions from the class $\Phi_{lsc}$ which posses the intersection property. 
\begin{proposition}
	\label{prop1n}
	Let $X$ be a Hilbert space, $f,g:X\rightarrow\bar{\mathbb{R}}$ be a proper $\Phi_{lsc}$-convex functions, $\alpha\in\mathbb{R}$ and $\varepsilon\geq0$. Assume that $\bar{x}\in \text{dom}(f)\cap \text{dom}(g)$ and $\bar{x}\in [f\geq \alpha]\cap[g\geq \alpha]$.
	
	If $f$ and $g$ satisfy the $ZS(\varepsilon, \bar{x})$ condition then, there exist
	$  \varphi_{1}\in\text{supp} (f), \ \varphi_{2}\in\text{supp} (g)$ for
	which the intersection property holds at the level  $ \alpha-\varepsilon$ (Definition \ref{ip}).
\end{proposition}
\begin{Proof}
	By Remark \ref{rem1}a) we only need to consider the case where $(a_1,v_1)\in \partial^\varepsilon_{lsc} f(\bar{x}) $ and $(a_2,v_2)\in\partial^\varepsilon_{lsc} g(\bar{x}) $ such that $\lambda v_1+\mu v_2=0$, $\lambda a_1+\mu a_2=0$ and
	$\lambda +\mu=1$ .
	
	If $\lambda =0$, then $0\in \partial^\varepsilon_{lsc} g(\bar{x})$ and we have $g(x)\geq g(\bar{x})-\varepsilon$ for all $x\in X$. By assumption $g(\bar{x})\geq \alpha$, so for $\varphi_1 \equiv \alpha-\varepsilon$, we have $\varphi_1\in \text{supp} (g)$, and $\varphi_1$ and any function $\varphi_2\in \text{supp}(f)$ have the intersection property at the level $\alpha-\varepsilon$, since $[\varphi_1<\alpha-\varepsilon]=\emptyset$.
	By the similar reasoning, we get the desired conclusion if $\mu=0$.
	
	Now assume that  $\lambda>0$ and $\mu> 0$, this implies that $a_1=a_2=0$, since $a_1,a_2\geq 0$. Let
	$$
	\varphi_1(x):=\langle v_1, x- \bar{x}\rangle +f(\bar{x})-\varepsilon\ \ \ \ \mbox{and}\ \ \ \
	\varphi_2(x):=\langle v_2, x- \bar{x}\rangle +g(\bar{x}) -\varepsilon
	$$ 
	for all $x\in X$. It is obvious that $\varphi_1\in \text{supp} (f)$ and $\varphi_2\in \text{supp} (g)$. Now we show that $\varphi_1$ and $\varphi_2$ have the intersection property at the level $\alpha-\varepsilon$. Let $x_1\in [\varphi_1<\alpha-\varepsilon]$, we have
	$$
	\begin{array}{l}
	\varphi_1(x_1)<\alpha-\varepsilon \ \ \Leftrightarrow \\
	\langle v_1, x_1- \bar{x}\rangle +f(\bar{x})-\varepsilon<\alpha -\varepsilon \ \  \Leftrightarrow \\
	\langle v_1, x_1- \bar{x}\rangle <\alpha-f(\bar{x})
	\end{array}
	$$
	By assumption that $\bar{x}\in [f\geq \alpha]\cap [g\geq \alpha] $, we have $ \alpha-f(\bar{x})\leq 0$, so
	$$
	\langle v_1, x_1- \bar{x}\rangle <0.
	$$
	Since, $x_1$ was chosen arbitrary, we get that  $ [\varphi_1<\alpha-\varepsilon]\subset [\langle v_1, \cdot- \bar{x}\rangle<0 ]$. By similar calculations we get  $ [\varphi_2<\alpha-\varepsilon]\subset [\langle v_2, \cdot- \bar{x}\rangle<0 ]$. Using the fact that $\lambda v_1+\mu v_2=0$ we have
	$$
	[\langle v_2, \cdot- \bar{x}\rangle<0 ]=[-\langle v_1, \cdot- \bar{x}\rangle<0 ]=[\langle v_1, \cdot- \bar{x}\rangle>0 ].
	$$
	
	Which means
	$$
	[\varphi_1<\alpha-\varepsilon] \cap [\varphi_2<\alpha-\varepsilon]=\emptyset.
	$$
\end{Proof}
\vskip 0.1 true in

The following Proposition shows that the intersection property (Definition \ref{ip}) of the functions in the support sets of $f$ and $g$ implies the zero subgradient  condition for $f$ and $g$. 
\begin{proposition}
	\label{nes}
	Let $f,g:X\rightarrow\bar{\mathbb{R}}$ be $\Phi_{lsc}$-convex functions, $\alpha\in\mathbb{R}$. 
	
	If there exist $  \varphi_{1}\in\text{supp} (f), \ \varphi_{2}\in\text{supp} (g)$ which  
	have the intersection property at the level $\alpha$, then for every $ \varepsilon>0$  there exist $x_{1} \in \textnormal{dom}(f)$, $x_{2}\in \textnormal{dom}(g)$ such that
	functions $f$ and $g$ satisfy the $ZS(\varepsilon,x_1,x_2)$ condition.
\end{proposition}
\begin{Proof}
	Let the intersection property (Definition \ref{ip}) holds for $\varphi_1\in\textnormal{supp}(f)$ and $\varphi_2\in\textnormal{supp}(g)$ at the level $\alpha$.
	Since  $\varphi_1\in\textnormal{supp}(f)$ and $\varphi_2\in\textnormal{supp}(g)$, we have
	$$
	\inf\limits_{x\in X}\{ f(x)-\varphi_{1}(x)\}=:d_1\ge 0 \ \ \text{and} \ \ \inf\limits_{x\in X} \{g(x)-\varphi_{2}(x)\}=:d_2\ge 0.
	$$
	Let $\varepsilon>0$, there exist  $x_{1} \in \textnormal{dom}(f)$ and $x_{2}\in \textnormal{dom}(g)$
	such that
	$$
	f(x_{1})-\varphi_{1}(x_{1})<d_{1}+\varepsilon \ \ \text{and} \ \ g(x_{2})-\varphi_{2}(x_{2})<d_{2}+\varepsilon
	$$
	and, we also have
	$$
	f(x)-\varphi_{1}(x)\geq d_{1} \ \ \text{and} \ \ g(x)-\varphi_{2}(x)\geq d_{2} \ \ \ \ \ \forall \ \ x\in X
	$$
	Since $\varphi_1,\varphi_2\in \Phi_{lsc}$, there exists $a_1,a_2\in \mathbb{R}_{+}$, $v_1,v_2\in X^*$ and $c_1,c_2$ such that $\varphi_1(\cdot)=-a_1\|\cdot\|^2+\langle v_1,\cdot \rangle +c_1$  and $\varphi_2(\cdot)=-a_2\|\cdot\|^2+\langle v_2,\cdot \rangle +c_2$. We have
	$$
	f(x)-f(x_1)\geq -a_1\|x\|^2+a_1\|x_1\|^2+\langle v_1,x-x_1 \rangle -\varepsilon \ \ \ \ \ \forall \ \ x\in X
	$$
	and 
	$$
	g(x)-g(x_2)\geq -a_2\|x\|^2+a_2\|x_1\|^2+\langle v_2,x-x_2 \rangle -\varepsilon \ \ \ \ \ \forall \ \ x\in X.
	$$
	Above inequalities mean that $(a_1,v_1)\in \partial^\varepsilon_{lsc} f(x_1)$ and $(a_2,v_2)\in \partial^\varepsilon_{lsc} g(x_2)$.
	
	Now we show that $
	0\in\mbox{\textnormal{co}}(\partial^{\varepsilon}_{lsc}f(x_1)\cup\partial^{\varepsilon}_{lsc} g(x_2)). 
	$ Since $\varphi_1,\varphi_2$ have the intersection property at the level $\alpha$, we consider the following cases:
	
	1. If $[\varphi_1<\alpha]=\emptyset$ then $a_1=0$, $v_1=0$, so $0\in \partial_{lsc}^{\varepsilon}f(x_1)$. 
	
	2. Analogously, if $[\varphi_2<\alpha]=\emptyset$, we get that $0\in \partial_{lsc}^{\varepsilon}f(x_2)$.
	
	3. Assume now, that $[\varphi_1<\alpha]\neq\emptyset$, $[\varphi_2<\alpha]\neq\emptyset$. If $a_1>0$ and $a_2>0$, then $\lim\limits_{\|x\|\rightarrow +\infty}\varphi_1(x)=\lim\limits_{\|x\|\rightarrow +\infty}\varphi_2(x)=-\infty$, this means that we can always find $x_0\in [\varphi_1<\alpha]\cap[\varphi_2<\alpha] $, which is a contradiction with the assumption that $\varphi_1$ and $\varphi_2$ have the intersection property at the level $\alpha$.
	
	If $a_1=0$ and $a_2>0$ we have $\lim\limits_{\|x\|\rightarrow +\infty}\varphi_2(x)=-\infty$. This means that, there exists $\delta>0$ such that for all  $x\in X$, $\|x\|>\delta$ we have
	$$
	\varphi_2(x)<\alpha
	$$
	Since $\varphi_2$ is affine function, in the set of all $x'$ such that
	$$
	\varphi_2(x')<\alpha
	$$
	we can always find an element $x'$ such that $\|x'\|>\delta$. This means that $ [\varphi_1<\alpha]\cap[\varphi_2<\alpha] \neq \emptyset$. A contradiction.
	
	If $a_1=a_2=0$, functions take the form $\varphi_1(\cdot)=\langle v_1,\cdot\rangle +c_1$ and $\varphi_2(\cdot)=\langle v_2,\cdot\rangle +c_2$. By assumption, $\varphi_1$ and $\varphi_2$ have the intersection property at the level $\alpha$ , so there exists $\lambda \in (0,1)$ such that
	$$
	\lambda \varphi_1(x)+(1-\lambda)\varphi_2(x)\geq \alpha \ \ \ \ \ \ \forall \ \ \  x\in X.
	$$
	This is equivalent to the following inequality
	$$
	\lambda v_1(x)+(1-\lambda)v_2(x)\geq \alpha -\lambda c_1-(1-\lambda)c_2\ \ \ \ \ \ \forall \ \ \  x\in X.
	$$ 
	Since $v_1,v_2$ are linear functionals this means that $\lambda v_1+(1-\lambda)v_2=0$.
\end{Proof}

The following example shows that we can not always have the equivalence between the zero subgradient condition and the intersection property for $\varepsilon=0$. 
\begin{example}
	Let $X=\mathbb{R}$ and $f,g:\mathbb{R}\rightarrow\mathbb{R}$ be such that 
	$$
	f(x)=\text{e}^x \ \ \ \ \ \ \ \ \ g(x)=-x^2+4.
	$$
	It is easy to see that $f$ and $g$ are $\Phi_{lsc}$-convex functions. Let $\alpha=0$. A function $\varphi_1(x)\equiv 0$ is in the support set of $f$, a function $\varphi_2(x)=-x^2$ is in the support set of $g$, and the intersection property holds for $\varphi_1$ and $\varphi_2$ at the level $\alpha=0$.
	It is easy to see that $0\notin	 \partial_{lsc}f(x_1)$ and $0\notin 	 \partial_{lsc}g(x_2)$ for any $x_1\in \text{dom}(f)$ and $x_2\in \text{dom}(g) $.
	We also have, that if $(a_2,v_2)\in  \partial_{lsc}g(x_2)$  then $a_2> 0$.
	This means that $0\notin \text{co}(\partial_{lsc}f(x_1) \cup  \partial_{lsc}g(x_2) )$ for any $x_1\in \text{dom}(f)$ and $x_2\in \text{dom}(g) $.
\end{example}

\section{Minimax theorem}
The main result of this section, the minimax theorem for $\Phi_{lsc}$-convex functions is based on the following general result proved in \cite{Syga2018} for an arbitrary set $X$ and class $\Phi$.
\begin{theorem}(Theorem 2.1 \cite{Syga2018}).
	\label{min-max}
	Let $X$ be a Hilbert space and $Y$ be a real vector space and let $a:X\times Y\rightarrow\bar{\mathbb{R}}$. Assume that for any $y\in Y$ the  function $a(\cdot,y):X\rightarrow\bar{\mathbb{R}}$ 
	is proper $\Phi_{lsc}$-convex on $X$ and for any $x\in X$ the function $a(x,\cdot):Y\rightarrow\bar{\mathbb{R}}$  is concave (in the classical sense)on $Y$. The following conditions are equivalent:
	\begin{description}
		\item [{\em (i)}] for every $\alpha\in\mathbb{R}$, $\alpha < \inf\limits_{x\in X} \sup\limits_{y\in Y} a(x,y)$, there exist $y_{1}, y_{2}\in Y$ and $\varphi_{1}\in \textnormal{supp } a(\cdot, y_{1})$, $\varphi_{2}\in \textnormal{supp } a(\cdot, y_{2})$ such that the intersection property holds for $\varphi_{1}$ and $\varphi_{2}$ on $X$ at the level $\alpha$,
		\item [{\em (ii)}] $\sup\limits_{y\in Y} \inf\limits_{x\in X} a(x,y)=\inf\limits_{x\in X} \sup\limits_{y\in Y} a(x,y).$
	\end{description}
\end{theorem}
In the above theorem the intersection property is sufficient and necessary condition for the minimax equality to hold. Taking into account the results from section \ref{grad} we show that our new condition, the zero subgradient condition, involving the $\Phi_{lsc}^{\varepsilon}$-subdifferential  is a sufficient condition for minimax equality to hold.
\begin{theorem}
	\label{new_min_max_1}
	Let $X,Y$ be a Hilbert spaces and  $a:X\times Y\rightarrow\bar{\mathbb{R}}$ be such that for any $y\in Y$ the  function $a(\cdot,y):X\rightarrow\mathbb{R}$ 
	is $\Phi_{lsc}$-convex on $X$ and for any $x\in X$ the function $a(x,\cdot):Y\rightarrow\mathbb{R}$  is concave on $Y$.

	If for every $\beta <\inf\limits_{x\in X} \sup\limits_{y\in Y} a(x,y)$ and every $\varepsilon>0$ there exist $y_{1}, y_{2}\in Y$ and $ \bar {x} \in [a(\cdot, y_{1})\geq \beta]\cap [a(\cdot, y_{2})\geq \beta]$ such that the functions $a(\cdot, y_{1})$ and $a(\cdot, y_{2})$ satisfy the $ZS(\varepsilon, \bar{x})$ condition then,
	$$\sup\limits_{y\in Y} \inf\limits_{x\in X} a(x,y)=\inf\limits_{x\in X} \sup\limits_{y\in Y} a(x,y).$$
\end{theorem}

\begin{Proof}
	We show that the condition $(i)$ of Theorem \ref{min-max} holds. Let $\alpha <\inf\limits_{x\in X} \sup\limits_{y\in Y} a(x,y)$. Let $\beta$ be such that $\alpha <\beta< \inf\limits_{x\in X} \sup\limits_{y\in Y} a(x,y)$ and $\varepsilon=\beta-\alpha>0$. By assumption, there exist $y_{1}, y_{2}\in Y$ and $ \bar {x} \in [a(\cdot, y_{1})\geq \beta]\cap [a(\cdot, y_{2})\geq \beta]$ such that
	$$
	0 \in \text{co}(\partial^{\varepsilon}_{lsc} a(\cdot, y_{1})(\bar{x})\cup \partial^{\varepsilon}_{lsc} a(\cdot, y_{2})(\bar{x})).
	$$
	From Proposition \ref{prop1n} we get that there exist $\varphi_{1}\in \textnormal{supp } a(\cdot, y_{1})$, $\varphi_{2}\in \textnormal{supp } a(\cdot, y_{2})$ such that the intersection property holds for $\varphi_{1}$ and $\varphi_{2}$ on $X$ at the level $\beta-\varepsilon$. Since, $\varepsilon=\beta-\alpha$, the intersection property holds at the level $\alpha$. 
\end{Proof}

Let us note that, if $\inf\limits_{x\in X} \sup\limits_{y\in Y} a(x,y)=-\infty$ then, the equality
$\sup\limits_{y\in Y} \inf\limits_{x\in X} a(x,y)=\inf\limits_{x\in X} \sup\limits_{y\in Y} a(x,y)$ always holds. If  we assume that $\inf\limits_{x\in X} \sup\limits_{y\in Y} a(x,y)<+\infty$ and $\beta = \inf\limits_{x\in X} \sup\limits_{y\in Y} a(x,y)$ we  get the following, stronger result.
\begin{theorem}
	\label{new_min_max_2}
	Let $X,Y$ be a Hilbert spaces and  $a:X\times Y\rightarrow\bar{\mathbb{R}}$ be   such that for any $y\in Y$ the  function $a(\cdot,y):X\rightarrow\mathbb{R}$ 
	is $\Phi_{lsc}$-convex on $X$ and for any $x\in X$ the function $a(x,\cdot):Y\rightarrow\mathbb{R}$  is concave on $Y$.

	If there exist $y_{1}, y_{2}\in Y$ and $\bar{x}\in \text{dom}(a(\cdot, y_{1}) )\cap \text{dom }(a(\cdot, y_{2}))$,  $ \bar {x} \in [a(\cdot, y_{1})\geq \beta]\cap [a(\cdot, y_{2})\geq \beta]$ such that the functions $a(\cdot, y_{1})$ and $a(\cdot, y_{2})$ satisfy the $ZS(0, \bar{x})$ condition, then
	$$\sup\limits_{y\in Y} \inf\limits_{x\in X} a(x,y)=\inf\limits_{x\in X} \sup\limits_{y\in Y} a(x,y).$$
\end{theorem}

\begin{Proof}
	The proof is an immediate consequence of Proposition \ref{prop1n} and Theorem \ref{min-max}.
\end{Proof}

Now, we can use Proposition \ref{nes} to provide the necessary condition, involving $\Phi_{lsc}$-subdifferential.
\begin{theorem}
	\label{ness_min_max}
	Let $X,Y$ be a Hilbert spaces and  $a:X\times Y\rightarrow\bar{\mathbb{R}}$ be a function such that for any $y\in Y$ the  function $a(\cdot,y):X\rightarrow\mathbb{R}$ 
	is $\Phi_{lsc}$-convex on $X$ and for any $x\in X$ the function $a(x,\cdot):Y\rightarrow\mathbb{R}$  is concave on $Y$. 
	
	If 	$$\sup\limits_{y\in Y} \inf\limits_{x\in X} a(x,y)=\inf\limits_{x\in X} \sup\limits_{y\in Y} a(x,y),$$
	then for every $\alpha <\inf\limits_{x\in X} \sup\limits_{y\in Y} a(x,y)$  there exist $y_{1}, y_{2}\in Y$ such that for every $\varepsilon>0$ there exist $x_{1} \in \textnormal{dom}(a(\cdot, y_{1}))$, $x_{2}\in \textnormal{dom}(a(\cdot, y_{2}))$ such that the functions  $a(\cdot, y_{1})$ and $a(\cdot, y_{2})$ satisfy the $ZS(\varepsilon,x_1,x_2)$.
\end{theorem}
\begin{Proof}
	The proof is an immediate consequence of Proposition \ref{nes} and Theorem \ref{min-max}.
\end{Proof}

It was established in the section 3 that even differentiable functions can have an empty $\Phi_{lsc}$-subdifferential, hence the zero subgradient condition for $\Phi_{lsc}$-convex functions seems to be strong one. We now formulate the minimax theorem for paraconvex functions, for which the nonemptiness  of  $\Phi_{lsc}$-subdifferential is assured.

\begin{theorem}
	\label{new_min_max_3}
	Let $X,Y$ be a Hilbert spaces and  $a:X\times Y\rightarrow\bar{\mathbb{R}}$ be   such that for any $y\in Y$ the  function $a(\cdot,y):X\rightarrow\mathbb{R}$ 
	is paraconvex on $X$ and for any $x\in X$ the function $a(x,\cdot):Y\rightarrow\mathbb{R}$  is concave on $Y$.

	If there exist $y_{1}, y_{2}\in Y$ and $\bar{x}\in \text{int dom}(a(\cdot, y_{1}) )\cap \text{int dom}(a(\cdot, y_{2}))$,  $ \bar {x} \in [a(\cdot, y_{1})\geq \beta]\cap [a(\cdot, y_{2})\geq \beta]$ such that the functions $a(\cdot, y_{1})$ and $a(\cdot, y_{2})$ satisfy the $ZS(0, \bar{x})$ condition, then
	$$\sup\limits_{y\in Y} \inf\limits_{x\in X} a(x,y)=\inf\limits_{x\in X} \sup\limits_{y\in Y} a(x,y).$$
\end{theorem}

\bibliographystyle{plain} 
\bibliography{bibn} 


\end{document}